# ASYMPTOTICS OF INPUT-CONSTRAINED BINARY SYMMETRIC CHANNEL CAPACITY

By Guangyue Han[1] and Brian Marcus

*University of Hong Kong and University of British Columbia*

We study the classical problem of noisy constrained capacity in the case of the binary symmetric channel (BSC), namely, the capacity of a BSC whose inputs are sequences chosen from a constrained set. Motivated by a result of Ordentlich and Weissman [In *Proceedings of IEEE Information Theory Workshop* (2004) 117–122], we derive an asymptotic formula (when the noise parameter is small) for the entropy rate of a hidden Markov chain, observed when a Markov chain passes through a BSC. Using this result, we establish an asymptotic formula for the capacity of a BSC with input process supported on an irreducible finite type constraint, as the noise parameter tends to zero.

**1. Introduction and background.** Let $X, Y$ be discrete random variables with alphabet $\mathcal{X}, \mathcal{Y}$ and joint probability mass function $p_{X,Y}(x,y) \stackrel{\triangle}{=} P(X = x, Y = y)$, $x \in \mathcal{X}, y \in \mathcal{Y}$ [for notational simplicity, we will write $p(x,y)$ rather than $p_{X,Y}(x,y)$, similarly $p(x), p(y)$ rather than $p_X(x), p_Y(y)$, resp., when it is clear from the context]. The *entropy* $H(X)$ of the discrete random variable $X$, which measures the level of uncertainty of $X$, is defined as (in this paper log is taken to mean the natural logarithm)

$$H(X) = -\sum_{x \in \mathcal{X}} p(x) \log p(x).$$

The *conditional entropy* $H(Y|X)$, which measures the level of uncertainty of $Y$ given $X$, is defined as

$$H(Y|X) = \sum_{x \in \mathcal{X}} p(x) H(Y|X = x)$$

Received September 2008.

[1]Supported by the University of Hong Kong under Grant No. 200709159007 and supported by the Research Grants Council of the Hong Kong Special Administrative Region, China under Grant No. HKU 701708P.

$$\stackrel{\triangle}{=} - \sum_{x \in \mathcal{X}} p(x) \sum_{y \in \mathcal{Y}} p(y|x) \log p(y|x) = - \sum_{x \in \mathcal{X}, y \in \mathcal{Y}} p(x,y) \log p(y|x).$$

The definitions above naturally include the case when $X, Y$ are vector-valued variables, for example, $X = X_k^\ell \stackrel{\triangle}{=} (X_k, X_{k+1}, \ldots, X_\ell)$, a sequence of discrete random variables.

For a left-infinite discrete stationary stochastic process $X = X_{-\infty}^0 \stackrel{\triangle}{=} \{X_i : i = 0, -1, -2, \ldots\}$, the *entropy rate* of $X$ is defined to be

$$(1.1) \qquad H(X) = \lim_{n \to \infty} \frac{1}{n+1} H(X_{-n}^0),$$

where $H(X_{-n}^0)$ denotes the entropy of the vector-valued random variable $X_{-n}^0$. Given another stationary process $Y = Y_{-\infty}^0$, we similarly define the conditional entropy rate

$$(1.2) \qquad H(Y|X) = \lim_{n \to \infty} \frac{1}{n+1} H(Y_{-n}^0 | X_{-n}^0).$$

A simple monotonicity argument in page 64 of [8] shows the existence of the limit in (1.1). Using the chain rule for entropy (see page 21 of [8]), we obtain

$$H(Y_{-n}^0 | X_{-n}^0) = H(X_{-n}^0, Y_{-n}^0) - H(X_{-n}^0),$$

and so we can apply the same argument to the processes $(X, Y)$ and $X$ to obtain the limit in (1.2).

If $Y = Y_{-\infty}^0$ is a stationary finite-state Markov chain, then $H(Y)$ has a simple analytic form. Specifically, denoting by $\Delta$ the transition probability matrix of $Y$, we have

$$(1.3) \qquad H(Y) = H(Y_0 | Y_{-1}) = - \sum_{i,j} P(Y_0 = i) \Delta(i,j) \log \Delta(i,j).$$

A function $Z = Z_{-\infty}^0$ of the stationary Markov chain $Y$ with the form $Z_i = \Phi(Y_i)$ is called a *hidden Markov chain*; here $\Phi$ is a function defined on the alphabet of $Y_i$, taking values in the alphabet of $Z_i$. We often write $Z = \Phi(Y)$. Hidden Markov chains are typically not Markov.

For a hidden Markov chain $Z$, the entropy rate $H(Z)$ was studied by Blackwell [6] as early as 1957, where the analysis suggested the intrinsic complexity of $H(Z)$ as a function of the process parameters. He gave an expression for $H(Z)$ in terms of a measure $Q$ on a simplex, obtained by solving an integral equation dependent on the parameters of the process. However, the measure is difficult to extract from the equation in any explicit way, and the entropy rate is difficult to compute.

Recently, the problem of computing the entropy rate of a hidden Markov chain has drawn much interest, and many approaches have been adopted to tackle this problem. These include asymptotic expansions as Markov chain parameters tend to extremes [14, 17, 18, 22, 23, 34, 35], analyticity results [13], variations on a classical bound [9] and efficient Monte Carlo methods [2, 27, 31]; and connections with the top Lyapunov exponent of a random matrix product have been observed [11, 15, 16, 17], relating to earlier work on Lyapunov exponents [4, 25, 26, 28].

Of particular interest are hidden Markov chains which arise as output processes of noisy channels. For example, the *binary symmetric channel with crossover probability* $\varepsilon$ [denoted BSC($\varepsilon$)] is an object which transforms input processes to output processes by means of a fixed i.i.d. binary noise process $E = \{E_n\}$ with $p_{E_n}(0) = 1 - \varepsilon$ and $p_{E_n}(1) = \varepsilon$. Specifically, given an arbitrary binary input process $X = \{X_n\}$, which is independent of $E$, define at time $n$,

$$Z_n(\varepsilon) = X_n \oplus E_n,$$

where $\oplus$ denotes binary addition modulo 2; then $Z_\varepsilon = \{Z_n(\varepsilon)\}$ is the output process corresponding to $X$.

When the input $X$ is a stationary Markov chain, the output $Z_\varepsilon$ can be viewed as a hidden Markov chain by appropriately augmenting the state space of $X$ [10]; specifically, in the case that $X$ is a first order binary Markov chain with transition probability matrix

$$\Pi = \begin{bmatrix} \pi_{00} & \pi_{01} \\ \pi_{10} & \pi_{11} \end{bmatrix},$$

then $Y_\varepsilon = \{Y_n(\varepsilon)\} = \{(X_n, E_n)\}$ is jointly Markov with transition probability matrix

$$\Delta = \begin{bmatrix} \begin{array}{c|cccc} y & (0,0) & (0,1) & (1,0) & (1,1) \\ \hline (0,0) & \pi_{00}(1-\varepsilon) & \pi_{00}\varepsilon & \pi_{01}(1-\varepsilon) & \pi_{01}\varepsilon \\ (0,1) & \pi_{00}(1-\varepsilon) & \pi_{00}\varepsilon & \pi_{01}(1-\varepsilon) & \pi_{01}\varepsilon \\ (1,0) & \pi_{10}(1-\varepsilon) & \pi_{10}\varepsilon & \pi_{11}(1-\varepsilon) & \pi_{11}\varepsilon \\ (1,1) & \pi_{10}(1-\varepsilon) & \pi_{10}\varepsilon & \pi_{11}(1-\varepsilon) & \pi_{11}\varepsilon \end{array} \end{bmatrix}$$

and $Z_\varepsilon = \{Z_n(\varepsilon)\}$ is a hidden Markov chain with $Z_n(\varepsilon) = \Phi(Y_n(\varepsilon))$, where $\Phi$ maps states $(0,0)$ and $(1,1)$ to 0 and maps states $(0,1)$ and $(1,0)$ to 1.

In Section 2 we give asymptotics for the entropy rate of a hidden Markov chain, obtained by passing a binary Markov chain, of arbitrary order, through BSC($\varepsilon$) as the noise $\varepsilon$ tends to zero. In Section 2.1 we review, from [18], the result when the transition probabilities are strictly positive. In Section 2.2 we develop the formula when some transition probabilities are zero (which is our main focus), thereby generalizing a specific result from [23].

The remainder of the paper is devoted to asymptotics for noisy constrained channel capacity. The *capacity* of the (unconstrained) BSC($\varepsilon$) is defined

(1.4) $$C(\varepsilon) = \lim_{n \to \infty} \sup_{X_{-n}^0} \frac{1}{n+1}(H(Z_{-n}^0(\varepsilon)) - H(Z_{-n}^0(\varepsilon)|X_{-n}^0));$$

here $X_{-n}^0$ is a finite-length input process from time $-n$ to 0 and $Z_{-n}^0(\varepsilon)$ is the corresponding output process. Seminal results of information theory, due to Shannon [30], include the following: (1) the capacity is the optimal rate of transmission possible with arbitrarily small probability of error, and (2) the capacity can be explicitly computed: $C(\varepsilon) = 1 - H(\varepsilon)$, where $H(\varepsilon)$ is the binary entropy function defined as

$$H(\varepsilon) = \varepsilon \log 1/\varepsilon + (1 - \varepsilon) \log 1/(1 - \varepsilon).$$

Generally speaking, it is very difficult to calculate the capacity of a generic channel. For a discrete memoryless channel without input-constraints, the Blahut–Arimoto algorithm [1, 7] can be applied to approximate the capacity numerically. A generalized Blahut–Arimoto algorithm has been proposed to numerically compute the local maximum mutual information rate of a finite state machine channel [32]. We are interested in *input-constrained* channel capacity, i.e., the capacity of $\mathrm{BSC}(\varepsilon)$, where the possible inputs are constrained, described as follows.

Let $\mathcal{X} = \{0, 1\}$, $\mathcal{X}^*$ denote all the finite length binary words, and $\mathcal{X}^n$ denote all the binary words with length $n$. A binary *finite type* constraint [20, 21] $\mathcal{S}$ is a subset of $\mathcal{X}^*$ defined by a finite set (denoted by $\mathcal{F}$) of forbidden words; in other words, any element in $\mathcal{S}$ does not contain any element in $\mathcal{F}$ as a contiguous subsequence. A prominent example is the $(d, k)$-RLL constraint $\mathcal{S}(d, k)$, which forbids any sequence with fewer than $d$ or more than $k$ consecutive zeros in between two 1's. For $\mathcal{S}(d, k)$ with $k < \infty$, a forbidden set $\mathcal{F}$ is:

$$\mathcal{F} = \{1\underbrace{0\cdots0}_{l}1 : 0 \le l < d\} \cup \{\underbrace{0\cdots0}_{k+1}\}.$$

When $k = \infty$, one can choose $\mathcal{F}$ to be

$$\mathcal{F} = \{1\underbrace{0\cdots0}_{l}1 : 0 \le l < d\};$$

in particular, when $d = 1, k = \infty$, $\mathcal{F}$ can be chosen to be $\{11\}$. These constraints on input sequences arise in magnetic recording in order to eliminate the most damaging error events [21].

We will use $\mathcal{S}_n$ to denote the subset of $\mathcal{S}$ consisting of words with length $n$. A finite type constraint $\mathcal{S}$ is *irreducible* if for any $u, v \in \mathcal{S}$, there is a $w \in \mathcal{S}$ such that $uwv \in \mathcal{S}$.

For a finite binary stochastic (not necessarily stationary) process $X = X_{-n}^0$, define the set of *allowed* words with respect to $X$ as

$$\mathcal{A}(X_{-n}^0) = \{w_{-n}^0 \in \mathcal{X}^{n+1} : P(X_{-n}^0 = w_{-n}^0) > 0\}.$$

For a left-infinite binary stochastic (again not necessarily stationary) process $X = X_{-\infty}^0$, define the set of *allowed* words with respect to $X$ as

$$\mathcal{A}(X) = \{w_{-m}^0 \in \mathcal{X}^* : m \ge 0, P(X_{-m}^0 = w_{-m}^0) > 0\}.$$

For a constrained $\mathrm{BSC}(\varepsilon)$ with input sequences in $\mathcal{S}$, the *noisy constrained capacity* $C(\mathcal{S}, \varepsilon)$ is defined as

$$C(\mathcal{S}, \varepsilon) = \lim_{n \to \infty} \sup_{\mathcal{A}(X_{-n}^0) \subseteq \mathcal{S}} \frac{1}{n+1}(H(Z_{-n}^0(\varepsilon)) - H(Z_{-n}^0(\varepsilon)|X_{-n}^0)),$$

where again $Z_{-n}^0(\varepsilon)$ is the output process corresponding to the input process $X_{-n}^0$. Let $\mathbb{P}$ (resp. $\mathbb{P}_n$) denote the set of all left-infinite (resp. length $n$) stationary processes over the alphabet $\mathcal{X}$. Using the approach in Section 12.4 of [12], one can show that

$$
\begin{aligned}
C(\mathcal{S}, \varepsilon) &= \lim_{n \to \infty} \sup_{X_{-n}^0 \in \mathbb{P}_{n+1}, \mathcal{A}(X_{-n}^0) \subseteq \mathcal{S}} \frac{1}{n+1}(H(Z_{-n}^0(\varepsilon)) - H(Z_{-n}^0(\varepsilon)|X_{-n}^0)) \\
&= \sup_{X \in \mathbb{P}, \ \mathcal{A}(X) \subseteq \mathcal{S}} H(Z_\varepsilon) - H(Z_\varepsilon|X),
\end{aligned}
$$

(1.5)

where $Z_{-n}^0(\varepsilon), Z_\varepsilon$ are the output process corresponding to the input processes $X_{-n}^0, X$, respectively.

In Section 3 we apply the results of Section 2 to derive an asymptotic formula for capacity of the input-constrained $\mathrm{BSC}(\varepsilon)$ (again as $\varepsilon$ tends to zero) for any irreducible finite type input constraint. In Section 4 we consider the special case of the $(d, k)$-RLL constraint, and compute the coefficients of the asymptotic formulas.

Regarding prior work on $C(\mathcal{S}, \varepsilon)$, the best results in the literature have been in the form of bounds and numerical simulations based on producing random (and, hopefully, typical) channel output sequences (see, e.g., [3, 29, 33] and references therein). These methods allow for fairly precise numerical approximations of the capacity for given constraints and channel parameters.

For a more detailed introduction to entropy, capacity and related concepts in information theory, we refer to standard textbooks such as [8, 12].

## 2. Asymptotics of entropy rate.

Consider a $\mathrm{BSC}(\varepsilon)$ and suppose the input is an $m$th order irreducible Markov chain $X$ defined by the transition probabilities $P(X_t = a_0 | X_{t-m}^{t-1} = a_{-m}^{-1})$, $a_{-m}^0 \in \mathcal{X}^{m+1}$, here again $\mathcal{X} = \{0, 1\}$, and the output hidden Markov chain will be denoted by $Z_\varepsilon$.

### 2.1. *When transition probabilities of $X$ are all positive.* This case is treated in [18]:

THEOREM 2.1 ([18], Theorem 3). *If $P(X_t = a_0 | X_{t-m}^{t-1} = a_{-m}^{-1}) > 0$ for all $a_{-m}^0 \in \mathcal{X}^{n+1}$, the entropy rate of $Z_\varepsilon$ for small $\varepsilon$ is*

$$(2.1) \qquad H(Z_\varepsilon) = H(X) + g(X)\varepsilon + O(\varepsilon^2),$$

*where, denoting by $\bar{z}_i$ the Boolean complement of $z_i$, and*

$$\check{z}^{2m+1} = z_1 \cdots z_m \bar{z}_{m+1} z_{m+2} \cdots z_{2m+1},$$

*we have*

$$(2.2) \qquad g(X) = \sum_{z_1^{2m+1} \in \mathcal{X}_{2m+1}} P_X(z_1^{2m+1}) \log \frac{P_X(z_1^{2m+1})}{P_X(\check{z}_1^{2m+1})}.$$

We remark that the expression here for $g(X)$ is a familiar quantity in information theory, known as the Kullback–Liebler divergence; specifically, $g(X)$ is the divergence between the two distributions $P_X(z_1^{2m+1})$ and $P_X(\check{z}_1^{2m+1})$.

In [18] a complete proof is given for first-order Markov chains, as well as the sketch for the generalization to higher order Markov chains. Alternatively, after appropriately enlarging the state space of $X$ to convert the $m$th order Markov chain to a first order Markov chain, we can use Theorem 1.1 of [13] to show $H(Z_\varepsilon)$ is analytic with respect to $\varepsilon$ at $\varepsilon = 0$, and Theorem 2.5 of [14] to show that all the derivatives of $H(Z_\varepsilon)$ at $\varepsilon = 0$ can be computed explicitly (in principle) without taking limits. Theorem 2.1 does this explicitly (in fact) for the first derivative.

### 2.2. When transition probabilities of $X$ are not necessarily all positive.

First consider the case when $X$ is a binary first order Markov chain with the transition probability matrix

$$(2.3) \qquad \begin{bmatrix} 1-p & p \\ 1 & 0 \end{bmatrix},$$

where $0 \le p \le 1$. This process generates sequences satisfying the $(d,k) = (1, \infty)$-RLL constraint, which simply means that the string 11 is forbidden. Sequences generated by the output process $Z_\varepsilon$, however, will generally not satisfy the constraint. The probability of the constraint-violating sequences at the output of the channel is polynomial in $\varepsilon$, which will generally contribute a term $O(\varepsilon \log \varepsilon)$ to the entropy rate $H(Z_\varepsilon)$ when $\varepsilon$ is small. This was already observed for the probability transition matrix (2.3) in [23], where it is shown that

$$(2.4) \qquad H(Z_\varepsilon) = H(X) + \frac{p(2-p)}{1+p} \varepsilon \log 1/\varepsilon + O(\varepsilon)$$

as $\varepsilon \to 0$.

In the following we shall generalize formulas (2.1) and (2.4) and derive a formula for entropy rate of any hidden Markov chain $Z_\varepsilon$, obtained when passing a Markov chain $X$ of any order $m$ through a BSC($\varepsilon$). We will apply the Birch bounds [5], for $n \geq m$, which yield

(2.5)
$$H(Z_0(\varepsilon)|Z_{-n+m}^{-1}(\varepsilon), X_{-n}^{-n+m-1}, E_{-n}^{-n+m-1})$$
$$\leq H(Z_\varepsilon) \leq H(Z_0(\varepsilon)|Z_{-n}^{-1}(\varepsilon)).$$

Note that the lower bound is really just

$$H(Z_0(\varepsilon)|Z_{-n+m}^{-1}(\varepsilon), X_{-n}^{-n+m-1}),$$

since $Z_{-n+m}^0(\varepsilon)$, if conditioned on $X_{-n}^{-n+m-1}$, is independent of $E_{-n}^{-n+m-1}$ .

LEMMA 2.2. *For a stationary input process $X_{-n}^0$ and the corresponding output process $Z_{-n}^0(\varepsilon)$ through BSC($\varepsilon$) and $0 \leq k \leq n$,*

$$H(Z_0(\varepsilon)|Z_{-n+k}^{-1}(\varepsilon), X_{-n}^{-n+k-1})$$
$$= H(X_0|X_{-n}^{-1}) + f_n^k(X_{-n}^0)\varepsilon \log(1/\varepsilon) + g_n^k(X_{-n}^0)\varepsilon + O(\varepsilon^2 \log \varepsilon),$$

*where $f_n^k(X_{-n}^0)$ and $g_n^k(X_{-n}^0)$ are given by (2.8) and (2.9) below, respectively.*

PROOF. In this proof $w = w_{-n}^{-1}$, where $w_{-j}$ is a single binary bit, and we let $v$ denote a single binary bit. And we use the notation for probability:

$$p_{XZ}(w) = P(X_{-n}^{-n+k-1} = w_{-n}^{-n+k-1}, Z_{-n+k}^{-1}(\varepsilon) = w_{-n+k}^{-1}),$$

$$p_{XZ}(wv) = P(X_{-n}^{-n+k-1} = w_{-n}^{-n+k-1}, Z_{-n+k}^{-1}(\varepsilon) = w_{-n+k}^{-1}, Z_0(\varepsilon) = v)$$

and

$$p_{XZ}(v|w) = P(Z_0(\varepsilon) = v|Z_{-n+k}^{-1}(\varepsilon) = w_{-n+k}^{-1}, X_{-n}^{-n+k-1} = w_{-n}^{-n+k-1}).$$

We remark that the definition of $p_{XZ}$ does depend on $\varepsilon$ and how we partition $w_{-n}^{-1}$ according to $k$, however, we keep the dependence implicit for notational simplicity.

We split $H(Z_0(\varepsilon)|Z_{-n+k}^{-1}(\varepsilon), X_{-n}^{-n+k-1})$ into five terms:

$$H(Z_0(\varepsilon)|Z_{-n+k}^{-1}(\varepsilon), X_{-n}^{-n+k-1})$$
$$= \sum_{wv \in \mathcal{A}(X)} -p_{XZ}(wv) \log(p_{XZ}(v|w))$$
$$+ \sum_{w \in \mathcal{A}(X), wv \notin \mathcal{A}(X)} -p_{XZ}(wv) \log(p_{XZ}(v|w))$$

(2.6)

$$+ \sum_{p_{XZ}(w)=\Theta(\varepsilon),\, p_{XZ}(wv)=\Theta(\varepsilon)} -p_{XZ}(wv)\log(p_{XZ}(v|w))$$

$$+ \sum_{p_{XZ}(w)=\Theta(\varepsilon),\, p_{XZ}(wv)=O(\varepsilon^2)} -p_{XZ}(wv)\log(p_{XZ}(v|w))$$

$$+ \sum_{p_{XZ}(w)=O(\varepsilon^2)} -p_{XZ}(wv)\log(p_{XZ}(v|w)),$$

here by $\alpha = \Theta(\beta)$, we mean, as usual, there exist positive constants $C_1, C_2$ such that $C_1|\beta| \le |\alpha| \le C_2|\beta|$, while by $\alpha = O(\beta)$, we mean there exists a positive constant $C$ such that $|\alpha| \le C|\beta|$; note that from

$$p_{XZ}(w)$$

$$= \sum_{u_{-n+k}^{-1}:\, w_{-n}^{-n+k-1}u_{-n+k}^{-1} \in \mathcal{A}(X_{-n}^{-1})} \Bigg( P(X_{-n}^{-n+k-1} = w_{-n}^{-n+k-1}, X_{-n+k}^{-1} = u_{-n+k}^{-1})$$

$$\times \prod_{j=-n+k}^{-1} p_E(u_j \oplus w_j) \Bigg),$$

we see that $p_{XZ}(w) = \Theta(\varepsilon)$ is equivalent to the statement that $w \notin \mathcal{A}(X_{-n}^{-1})$, and by flipping exactly one of the bits in $w_{-n+k}^{-1}$, one obtains, from $w$, a sequence in $\mathcal{A}(X_{-n}^{-1})$.

For the fourth term, we have

$$\sum_{p_{XZ}(w)=\Theta(\varepsilon),\, p_{XZ}(wv)=O(\varepsilon^2)} -p_{XZ}(wv)\log(p_{XZ}(v|w)) = O(\varepsilon^2 \log \varepsilon).$$

For the fifth term, we have

$$\sum_{p_{XZ}(w)=O(\varepsilon^2)} -p_{XZ}(wv)\log(p_{XZ}(v|w))$$

$$= \sum_{p_{XZ}(w)=O(\varepsilon^2)} -p_{XZ}(w) \sum_v p_{XZ}(v|w)\log(p_{XZ}(v|w))$$

$$\le (\log 2) \sum_{p_{XZ}(w)=O(\varepsilon^2)} p_{XZ}(w) = O(\varepsilon^2),$$

where we use the fact that $-\sum_v p_{XZ}(v|w)\log(p_{XZ}(v|w)) \le \log 2$ for any $w$. We conclude that the sum of the fourth term and the fifth term is $O(\varepsilon^2 \log \varepsilon)$.

For a binary sequence $u_{-n}^{-1}$, define $h_n^k(u_{-n}^{-1})$ to be

$$(2.7) \qquad h_n^k(u_{-n}^{-1}) = \sum_{j=1}^{n-k} p_X(u_{-n}^{-j-1} \bar{u}_{-j} u_{-j+1}^{-1}) - (n-k)p_X(u_{-n}^{-1}).$$

Note that with this notation, $h_n^k(w)$ and $h_{n+1}^k(wv)$ can be expressed as derivatives with respect to $\varepsilon$ at $\varepsilon = 0$:

$$h_n^k(w) = p_{XZ}'(w)|_{\varepsilon=0}, \qquad h_{n+1}^k(wv) = p_{XZ}'(wv)|_{\varepsilon=0}.$$

Then for the first term, we have

$$\sum_{wv \in \mathcal{A}(X)} -p_{XZ}(wv) \log(p_{XZ}(v|w))$$

$$= -\sum_{wv \in \mathcal{A}(X)} (p_X(wv) + h_{n+1}^k(wv)\varepsilon + O(\varepsilon^2))$$

$$\times \log\Big(p_X(v|w) + \frac{h_{n+1}^k(wv)p_X(w) - h_n^k(w)p_X(wv)}{p_X^2(w)}\varepsilon$$

$$+ O(\varepsilon^2)\Big)$$

$$= H(X_0|X_{-n}^{-1}) - \sum_{wv \in \mathcal{A}(X)} \Big(h_{n+1}^k(wv) \log p_X(v|w)$$

$$+ \frac{h_{n+1}^k(wv)p_X(w) - h_n^k(w)p_X(wv)}{p_X(w)}\Big)\varepsilon$$

$$+ O(\varepsilon^2).$$

For the second term, it is easy to check that for $w \in \mathcal{A}(X)$ and $wv \notin \mathcal{A}(X)$, $p_{XZ}(v|w) = \Theta(\varepsilon)$ and so

$$p_{XZ}(wv) = h_{n+1}^k(wv)\varepsilon + O(\varepsilon^2);$$

we then obtain

$$\sum_{w \in \mathcal{A}(X), wv \notin \mathcal{A}(X)} -p_{XZ}(wv) \log(p_{XZ}(v|w))$$

$$= -\sum_{w \in \mathcal{A}(X), wv \notin \mathcal{A}(X)} h_{n+1}^k(wv)\varepsilon \log \frac{h_{n+1}^k(wv)\varepsilon + O(\varepsilon^2)}{p_X(w)} + O(\varepsilon^2) \log \Theta(\varepsilon)$$

$$= \sum_{w \in \mathcal{A}(X), wv \notin \mathcal{A}(X)} h_{n+1}^k(wv)\varepsilon \log(1/\varepsilon)$$

$$- \Big(\sum_{w \in \mathcal{A}(X), wv \notin \mathcal{A}(X)} h_{n+1}^k(wv) \log \frac{h_{n+1}^k(wv)}{p_X(w)}\Big)\varepsilon + O(\varepsilon^2 \log \varepsilon).$$

For the third term, we have

$$\sum_{p_{XZ}(w)=\Theta(\varepsilon), p_{XZ}(wv)=\Theta(\varepsilon)} -p_{XZ}(wv) \log(p_{XZ}(v|w))$$

$$= - \sum_{p_{XZ}(w) = \Theta(\varepsilon), p_{XZ}(wv) = \Theta(\varepsilon)} (h_{n+1}^k(wv)\varepsilon + O(\varepsilon^2))$$

$$\times \log\left(\frac{h_{n+1}^k(wv)}{h_n^k(w)} + O(\varepsilon)\right)$$

$$= -\left(\sum_{p_{XZ}(w) = \Theta(\varepsilon), p_{XZ}(wv) = \Theta(\varepsilon)} h_{n+1}^k(wv) \log\left(\frac{h_{n+1}^k(wv)}{h_n^k(w)}\right)\right)\varepsilon + O(\varepsilon^2).$$

In summary, $H(Z_0(\varepsilon)|Z_{-n+k}^{-1}(\varepsilon), X_{-n}^{-n+k-1})$ can be rewritten as

$$H(Z_0(\varepsilon)|Z_{-n+k}^{-1}(\varepsilon), X_{-n}^{-n+k-1})$$

$$= H(X_0|X_{-n}^{-1}) + f_n^k(X_{-n}^0)\varepsilon \log(1/\varepsilon) + g_n^k(X_{-n}^0)\varepsilon + O(\varepsilon^2 \log \varepsilon),$$

where [see (2.7) for the definition of $h_n^k(\cdot)$]

$$(2.8) \quad \begin{aligned} f_n^k(X_{-n}^0) &= \sum_{w \in \mathcal{A}(X), wv \notin \mathcal{A}(X)} h_{n+1}^k(wv) \\ &= \sum_{w \in \mathcal{A}(X), wv \notin \mathcal{A}(X)} \left(\sum_{j=1}^{n-k} p_X(w_{-n}^{-j-1} \bar{w}_{-j} w_{-j+1}^{-1} v) + p_X(w_{-n}^{-1})\right), \end{aligned}$$

and

$$(2.9) \quad \begin{aligned} g_n^k(X_{-n}^0) = &- \sum_{wv \in \mathcal{A}(X)} \left(h_{n+1}^k(wv) \log p_X(v|w)\right. \\ &\qquad\qquad \left. + \frac{h_{n+1}^k(wv) p_X(w) - h_n^k(w) p_X(wv)}{p_X(w)}\right) \\ &- \sum_{w \in \mathcal{A}(X), wv \notin \mathcal{A}(X)} h_{n+1}^k(wv) \log \frac{h_{n+1}^k(wv)}{p_X(w)} \\ &- \sum_{p_{XZ}(w) = \Theta(\varepsilon), p_{XZ}(wv) = \Theta(\varepsilon)} h_{n+1}^k(wv) \log\left(\frac{h_{n+1}^k(wv)}{h_n^k(w)}\right). \qquad \square \end{aligned}$$

REMARK 2.3.   For any $\delta > 0$ and fixed $n$, the constant in $O(\varepsilon^2 \log \varepsilon)$ in Lemma 2.2 can be chosen uniformly on $\mathbb{P}_{n+1,\delta}$, where $\mathbb{P}_{n+1,\delta}$ denotes the set of binary stationary processes $X = X_{-n}^0$, such that, for all $w_{-n}^0 \in \mathcal{A}(X)$, we have $p_X(w) \geq \delta$.

THEOREM 2.4.   For an $m$th order Markov chain $X$ passing through a $BSC(\varepsilon)$, with $Z_\varepsilon$ as the output hidden Markov chain,

$$H(Z_\varepsilon) = H(X) + f(X)\varepsilon \log(1/\varepsilon) + g(X)\varepsilon + O(\varepsilon^2 \log \varepsilon),$$

where $f(X) = f_{2m}^0(X_{-2m}^0) = f_{2m}^m(X_{-2m}^0)$ and $g(X) = g_{3m}^0(X_{-3m}^0) = g_{3m}^m(X_{-3m}^0)$.

PROOF. We apply Lemma 2.2 to the Birch upper and lower bounds [equation (2.5)] of $H(Z_\varepsilon)$. For the upper bound, $k = 0$, we have, for all $n$,

$$H(Z_0(\varepsilon)|Z_{-n}^{-1}(\varepsilon)) = H(X_0|X_{-n}^{-1}) + f_n^0(X_{-n}^0)\varepsilon \log(1/\varepsilon)$$
$$+ g_n^0(X_{-n}^0)\varepsilon + O(\varepsilon^2 \log \varepsilon).$$

And for the lower bound, $k = m$, we have, for $n \geq m$,

$$H(Z_0(\varepsilon)|Z_{-n+m}^{-1}(\varepsilon), X_{-n}^{-n+m-1})$$
$$= H(X_0|X_{-n}^{-1}) + f_n^m(X_{-n}^0)\varepsilon \log(1/\varepsilon) + g_n^m(X_{-n}^0)\varepsilon + O(\varepsilon^2 \log \varepsilon).$$

The first term always coincides for the upper and lower bounds. When $n \geq m$, since $X$ is an $m$th order Markov chain,

$$H(X_0|X_{-n}^{-1}) = H(X_0|X_{-m}^{-1}) = H(X).$$

Let $w = w_{-n}^{-1}$, where $w_{-j}$ is a single bit, and $v$ denotes a single bit. If $w \in \mathcal{A}(X)$ and $wv \notin \mathcal{A}(X)$, then $p_X(w_{-m}^{-1}v) = 0$. It then follows that for an $m$th order Markov chain, when $n \geq 2m$,

$$(2.10) \qquad f_n^m(X_{-n}^0) = f_n^0(X_{-n}^0) = f_{2m}^0(X_{-2m}^0) = f_{2m}^m(X_{-2m}^0).$$

Now consider $g_n^k(X_{-n}^0)$. When $0 \leq k \leq m$, we have the following facts [for a detailed derivation of (2.11)–(2.13), see the Appendix]:

$$(2.11) \qquad \text{if } wv \in \mathcal{A}(X), \qquad p_X(v|w) = p_X(v|w_{-m}^{-1}) \qquad \text{for } n \geq m,$$

$$(2.12) \qquad \begin{array}{l} \text{if } w \in \mathcal{A}(X),\, wv \notin \mathcal{A}(X), \\[4pt] \dfrac{h_{n+1}^k(wv)}{p_X(w)} \text{ is constant (as function of } n \text{ and } k\text{) for } n \geq 2m, 0 \leq k \leq m, \end{array}$$

$$(2.13) \qquad \begin{array}{l} \text{if } p_{XZ}(w) = \Theta(\varepsilon), p_{XZ}(wv) = \Theta(\varepsilon), \\[4pt] \dfrac{h_{n+1}^k(wv)}{h_n^k(w)} \text{ is constant for } n \geq 3m, 0 \leq k \leq m. \end{array}$$

It then follows [see the derivations of (2.14)–(2.16) in the Appendix] that

$$(2.14) \qquad \sum_{wv \in \mathcal{A}(X)} \frac{h_{n+1}^k(wv)p_X(w) - h_n^k(w)p_X(wv)}{p_X(w)}$$

is constant (as a function of $n$) for $n \geq 2m, 0 \leq k \leq m$,

$$(2.15) \qquad \sum_{w \in \mathcal{A}(X), wv \notin \mathcal{A}(X)} h_{n+1}^k(wv) \log \frac{h_{n+1}^k(wv)}{p_X(w)}$$

is constant for $n \geq 2m, 0 \leq k \leq m$, and

$$(2.16) \quad \sum_{wv \in \mathcal{A}(X)} h_{n+1}^k(wv) \log p_X(v|w)$$

$$+ \sum_{p_{XZ}(w) = \Theta(\varepsilon), p_{XZ}(wv) = \Theta(\varepsilon)} h_{n+1}^k(wv) \log \frac{h_{n+1}^k(wv)}{h_n^k(w)}$$

is constant for $n \geq 3m, 0 \leq k \leq m$.

Consequently, we have

$$(2.17) \quad g_n^m(X_{-n}^0) = g_n^0(X_{-n}^0) = g_{3m}^0(X_{-3m}^0) = g_{3m}^m(X_{-3m}^0).$$

Let $f(X) = f_{2m}^0(X_{-2m}^0)$ and $g(X) = g_{3m}^0(X_{-3m}^0)$, then the theorem follows.
$\square$

REMARK 2.5. Note that this result applies in particular to the case when the transition probabilities of $X$ are all positive; thus, in this case the formula should reduce to that of Theorem 2.1. Indeed, when all transition probabilities of $X$ are positive, $f(X)$ vanishes since the summation in (2.8) is taken over an empty set; on the other hand, again from (2.8), if some of the transition probabilities of $X$ are zero, then $f(X)$ does not vanish [to see this, note that when $w \in \mathcal{A}(X), wv \notin \mathcal{A}(X)$, necessarily we will have $w\bar{v} \in \mathcal{A}(X)$]. The agreement of $g(X)$ with expression in Theorem 2.1 is a straightforward, but tedious, computation.

REMARK 2.6. Together with Remark 2.3, the proof of Theorem 2.4 implies that for any $\delta > 0$ and fixed $m$, the constant in $O(\varepsilon^2 \log \varepsilon)$ in Theorem 2.4 can be chosen uniformly on $\mathbb{Q}_{m,\delta}$, where $\mathbb{Q}_{m,\delta}$ denotes the set of all $m$th order Markov chains $X$ such that, whenever $w = w_{-m}^0 \in \mathcal{A}(X)$, we have $p_X(w) \geq \delta$.

REMARK 2.7. The error term in the formula of Theorem 2.4 cannot be improved, in the sense that, in some cases, the error term is dominated by a strictly positive constant times $\varepsilon^2 \log \varepsilon$.

As we showed in Theorem 2.4, the Birch upper bound with $n = 3m$ yields

$$H(Z_0(\varepsilon)|Z_{-n}^{-1}(\varepsilon)) = H(X) + f(X)\varepsilon \log(1/\varepsilon) + g(X)\varepsilon + O(\varepsilon^2 \log \varepsilon).$$

Together with (2.6), one checks that the $\Theta(\varepsilon^2 \log \varepsilon)$ term in the error term $O(\varepsilon^2 \log \varepsilon)$ is contributed by [see the second term in (2.6) with $k = 0$]

$$\sum_{w \in \mathcal{A}(X), wv \notin \mathcal{A}(X)} -p_Z(wv) \log(p_Z(v|w))$$

and [see the fourth term in (2.6) with $k = 0$]

$$\sum_{p_Z(w) = \Theta(\varepsilon), p_Z(wv) = O(\varepsilon^2)} -p_Z(wv) \log(p_Z(v|w)),$$

and this $\Theta(\varepsilon^2 \log \varepsilon)$ term does not vanish at least for certain cases. For instance, consider the input Markov chain $X$ with the following transition probability matrix:

$$\begin{bmatrix} 1-p & p \\ 1 & 0 \end{bmatrix},$$

where $0 < p < 1$. Then one checks that for this case, $m = 1, n = 3$, and the coefficient of the above-mentioned $\Theta(\varepsilon^2 \log \varepsilon)$ term takes the form of

$$\frac{1 - 6p + 7p^2 - p^3}{1 + p},$$

which is strictly positive for $p$ is close to 0.

**3. Asymptotics of capacity.** Consider a binary irreducible finite type constraint $\mathcal{S}$ defined by $\mathcal{F}$, which consists of forbidden words with length $\hat{m} + 1$. In general, there are many such $\mathcal{F}$'s corresponding to the same $\mathcal{S}$ with different lengths; here we may choose $\mathcal{F}$ to be the one with the smallest length $\hat{m} + 1$. And $\hat{m} = \hat{m}(\mathcal{S})$ is defined to be the *topological order* of the constraint $\mathcal{S}$. For example, the order of $\mathcal{S}(d, k)$, discussed in the introduction, is $k$ [20]. The topological order of a finite type constraint is analogous to the order of a Markov chain.

Recall from (1.5) that for an input-constrained BSC($\varepsilon$) with input sequences $X_{-n}^0$ in $\mathcal{S}$ and with the corresponding output $Z_{-n}^0(\varepsilon)$, the capacity can be written as

$$C(\mathcal{S}, \varepsilon) = \lim_{n \to \infty} \sup_{X_{-n}^0 \in \mathbb{P}_{n+1}, \mathcal{A}(X_{-n}^0) \subseteq \mathcal{S}} (1/(n+1)(H(Z_{-n}^0(\varepsilon)) - H(Z_{-n}^0(\varepsilon)|X_{-n}^0)))$$

Since the noise distribution is symmetric and the noise process $E$ is i.i.d. and independent of $X$, this can be simplified to

$$C(\mathcal{S}, \varepsilon) = \lim_{n \to \infty} \sup_{X_{-n}^0 \in \mathbb{P}_{n+1}, \mathcal{A}(X_{-n}^0) \subseteq \mathcal{S}} H(Z_{-n}^0(\varepsilon))/(n+1) - H(\varepsilon),$$

which can be rewritten as

$$C(\mathcal{S}, \varepsilon) = \lim_{n \to \infty} \sup_{X_{-n}^0 \in \mathbb{P}_{n+1}, \mathcal{A}(X_{-n}^0) \subseteq \mathcal{S}} H(Z_0(\varepsilon)|Z_{-n}^{-1}(\varepsilon)) - H(\varepsilon),$$

where we used the chain rule for entropy (see page 21 of [8])

$$H(Z_{-n}^0(\varepsilon)) = \sum_{j=0}^{n} H(Z_0(\varepsilon)|Z_{-j}^{-1}(\varepsilon)),$$

and the fact that (further) conditioning reduces entropy (see page 27 of [8])

$$H(Z_0(\varepsilon)|Z_{-j_1}^{-1}(\varepsilon)) \geq H(Z_0(\varepsilon)|Z_{-j_2}^{-1}(\varepsilon)) \qquad \text{for } j_1 \leq j_2.$$

Recall from (1.5) that

$$C(\mathcal{S}, \varepsilon) = \sup_{X \in \mathbb{P}, \mathcal{A}(X) \subseteq \mathcal{S}} H(Z_\varepsilon) - H(Z_\varepsilon | X).$$

Now let

$$\mathcal{H}_n(\mathcal{S}, \varepsilon) = \sup_{X_{-n}^0 \in \mathbb{P}_{n+1}, \mathcal{A}(X_{-n}^0) \subseteq \mathcal{S}} H(Z_0(\varepsilon)|Z_{-n}^{-1}(\varepsilon))$$

and

$$h_m(\mathcal{S}, \varepsilon) = \sup_{X \in \mathcal{M}_m, \mathcal{A}(X) \subseteq \mathcal{S}} H(Z_\varepsilon),$$

where $\mathcal{M}_m$ denotes the set of all $m$th order binary irreducible Markov chains; we then have the bounds for $C(\mathcal{S}, \varepsilon)$:

$$(3.1) \qquad h_m(\mathcal{S}, \varepsilon) - H(\varepsilon) \leq C(\mathcal{S}, \varepsilon) \leq \mathcal{H}_n(\mathcal{S}, \varepsilon) - H(\varepsilon).$$

Noting that

$$\sup_{X_{-n}^0 \in \mathbb{P}_{n+1}, \mathcal{A}(X_{-n}^0) \subsetneq \mathcal{S}_{n+1}} H(X_0|X_{-n}^{-1}) < \sup_{X_{-n}^0 \in \mathbb{P}_{n+1}, \mathcal{A}(X_{-n}^0) = \mathcal{S}_{n+1}} H(X_0|X_{-n}^{-1})$$

(here $\subsetneq$ means "proper subset of"), and $H(Z_0(\varepsilon)|Z_{-n}^{-1}(\varepsilon))$ are continuous at $\varepsilon = 0$, we conclude that, for $\varepsilon$ sufficiently small ($\varepsilon < \varepsilon_0$), one may choose $\delta > 0$ (here, $\delta$ depends on $n$ and $m$) such that

$$\mathcal{H}_n(\mathcal{S}, \varepsilon) = \sup_{X_{-n}^0 \in \mathbb{P}_{n+1,\delta}, \mathcal{A}(X_{-n}^0) = \mathcal{S}_{n+1}} H(Z_0(\varepsilon)|Z_{-n}^{-1}(\varepsilon)).$$

So from now on we only consider stationary processes $X = X_{-n}^0$ with $\mathcal{A}(X_{-n}^0) = \mathcal{S}_{n+1}$.

Now for a stationary process $X = X_{-n}^0$, define $\vec{p}_n$ as the following probability vector indexed by all the elements in $\mathcal{S}_{n+1}$:

$$\vec{p}_n = \vec{p}_n(X_{-n}^0) = (P(X_{-n}^0 = w_{-n}^0) : w_{-n}^0 \in \mathcal{S}_{n+1}).$$

To emphasize the dependence of $X_{-n}^0$ on $\vec{p}_n$, in the following, we shall rewrite $X_{-n}^0$ as $X_{-n}^0(\vec{p}_n)$. For an $m$th order binary irreducible Markov chain $X = X_{-\infty}^0$, slightly abusing the notation, define $\vec{p}_m$ as the following probability vector indexed by all the elements in $\mathcal{S}_{m+1}$,

$$\vec{p}_m = \vec{p}_m(X_{-\infty}^0) = (P(X_{-m}^0 = w_{-m}^0) : w_{-m}^0 \in \mathcal{S}_{m+1}).$$

Similarly, to emphasize the dependence of $X = X_{-\infty}^0$ on $\vec{p}_m$, in the following, we shall rewrite $X$ as $X_{\vec{p}_m}$. And we shall use $Z_{-n}^0(\vec{p}_n, \varepsilon)$ to denote the output process obtained by passing $X_{-n}^0(\vec{p}_n)$ through BSC($\varepsilon$), and use $Z_{\vec{p}_m, \varepsilon}$ to denote the output process obtained by passing $X_{\vec{p}_m}$ through BSC($\varepsilon$).

LEMMA 3.1.   *For any stationary process $X^0_{-n}(\vec{p}_n)$ with $\mathcal{A}(X^0_{-n}(\vec{p}_n)) = \mathcal{S}_{n+1}$, $H(X_0(\vec{p}_n)|X^{-1}_{-n}(\vec{p}_n))$, as a function of $\vec{p}_n$, has a negative definite Hessian matrix.*

PROOF.   Note that

$$H(X_0(\vec{p}_n)|X^{-1}_{-n}(\vec{p}_n)) = -\sum_{x^0_{-n} \in \mathcal{S}} p(x^0_{-n}) \log p(x_0|x^{-1}_{-n}).$$

For two different probability vectors $\vec{p}_n$ and $\vec{q}_n$, consider the convex combination

$$\vec{r}_n(t) = t\vec{p}_n + (1-t)\vec{q}_n,$$

where $0 \le t \le 1$. It suffices to prove that $H(X_0(\vec{r}_n(t))|X^{-1}_{-n}(\vec{r}_n(t)))$ has a strictly negative second derivative with respect to $t$. Now consider a single term in $H(X_0(\vec{p}_n)|X^{-1}_{-n}(\vec{p}_n))$:

$$-(t\vec{p}_n(x^0_{-n}) + (1-t)\vec{q}_n(x^0_{-n})) \log \frac{t\vec{p}_n(x^0_{-n}) + (1-t)\vec{q}_n(x^0_{-n})}{t\vec{p}_n(x^{-1}_{-n}) + (1-t)\vec{q}_n(x^{-1}_{-n})}.$$

Note that for two formal symbols $\alpha$ and $\beta$, if we assume $\alpha'' = 0$ and $\beta'' = 0$, the second order formal derivative of $\alpha \log \frac{\alpha}{\beta}$ can be computed as

$$\left(\alpha \log \frac{\alpha}{\beta}\right)'' = \left(\frac{\alpha'}{\sqrt{\alpha}} - \sqrt{\alpha}\frac{\beta'}{\beta}\right)^2.$$

It then follows that the second derivative of this term (with respect to $t$) can be calculated as

$$-\bigg(\frac{\vec{p}_n(x^0_{-n}) - \vec{q}_n(x^0_{-n})}{\sqrt{t\vec{p}_n(x^0_{-n}) + (1-t)\vec{q}_n(x^0_{-n})}}$$

$$- \sqrt{t\vec{p}_n(x^0_{-n}) + (1-t)\vec{q}_n(x^0_{-n})}\,\frac{\vec{p}_n(x^0_{-(n-1)}) - \vec{q}_n(x^0_{-(n-1)})}{t\vec{p}_n(x^0_{-(n-1)}) + (1-t)\vec{q}_n(x^0_{-(n-1)})}\bigg)^2.$$

That is, the expression above is always nonpositive, and is equal to 0 only if

$$\frac{\vec{p}_n(x^0_{-n}) - \vec{q}_n(x^0_{-n})}{t\vec{p}_n(x^0_{-n}) + (1-t)\vec{q}_n(x^0_{-n})} = \frac{\vec{p}_n(x^0_{-(n-1)}) - \vec{q}_n(x^0_{-(n-1)})}{t\vec{p}_n(x^0_{-(n-1)}) + (1-t)\vec{q}_n(x^0_{-(n-1)})},$$

which is equivalent to

$$(3.2) \qquad \begin{aligned} &P(X_0(\vec{p}_n) = x_0|X^{-1}_{-n}(\vec{p}_n) = x^{-1}_{-n}) \\ &= P(X_0(\vec{q}_n) = x_0|X^{-1}_{-n}(\vec{q}_n) = x^{-1}_{-n}). \end{aligned}$$

Since $\mathcal{S}$ is an irreducible finite type constraint and $\mathcal{A}(X^0_{-n}(\vec{p}_n)) = \mathcal{A}(X^0_{-n}(\vec{q}_n)) = \mathcal{S}_{n+1}$, the expression (3.2) cannot be true for every $x^0_{-n}$ unless $\vec{p}_n = \vec{q}_n$. So we conclude that the second derivative of $H(X_0(\vec{r}_n(t))|X_{-n}^{-1}(\vec{r}_n(t)))$ (with respect to $t$) is strictly negative. Thus, $H(X_0(\vec{p}_n)|X_{-n}^{-1}(\vec{p}_n))$, as a function of $\vec{p}_n$, has a strictly negative definite Hessian. $\quad\square$

For $m \geq \hat{m}$, over all $m$th order Markov chains $X_{\vec{p}_m}$ with $\mathcal{A}(X_{\vec{p}_m}) = \mathcal{S}$, $H(X_{\vec{p}_m})$ is maximized at some unique Markov chain $X_{\vec{p}_m^{\max}}$ (see [20], [24]). Moreover, $X_{\vec{p}^{\max}}$ does not depend on $m$ and is an $\hat{m}$th order Markov chain, so we will drop the subscript $m$ and use $X_{\vec{p}^{\max}}$ instead to denote $X_{\vec{p}_m^{\max}}$ for any $m \geq \hat{m}$. The same idea shows that over all stationary distributions $X_{-n}^0(\vec{p}_n)$ $(n \geq \hat{m})$ with $\mathcal{A}(X_{-n}^0(\vec{p}_n)) = \mathcal{S}_{n+1}$, $H(X_0(\vec{p}_n)|X_{-n}^{-1}(\vec{p}_n))$ is maximized at $\vec{p}_n^{\max}$, which corresponds to the above unique $X_{\vec{p}^{\max}}$ as well.

Note that $C(\mathcal{S}) = C(\mathcal{S}, 0)$ is equal to the *noiseless capacity* of the constraint $\mathcal{S}$. This quantity has been extensively studied, and several interpretations and methods for its explicit derivation are known (see, e.g., [21] and the extensive bibliography therein). It is well known that $C(\mathcal{S}) = H(X_{\vec{p}^{\max}})$ (see [20], [24]).

THEOREM 3.2.    *1. If $n \geq 3\hat{m}(\mathcal{S})$,*

$$\mathcal{H}_n(\mathcal{S}, \varepsilon) = C(\mathcal{S}) + f(X_{\vec{p}^{\max}})\varepsilon \log(1/\varepsilon) + g(X_{\vec{p}^{\max}})\varepsilon + O(\varepsilon^2 \log^2 \varepsilon).$$

*2. If $m \geq \hat{m}(\mathcal{S})$,*

$$h_m(\mathcal{S}, \varepsilon) = C(\mathcal{S}) + f(X_{\vec{p}^{\max}})\varepsilon \log(1/\varepsilon) + g(X_{\vec{p}^{\max}})\varepsilon + O(\varepsilon^2 \log^2 \varepsilon).$$

*Here, as defined in Theorem 2.4, $f(X_{\vec{p}^{\max}}) = f_{2\hat{m}}^0(X_{-2\hat{m}}^0(\vec{p}^{\max}))$ and $g(X_{\vec{p}^{\max}}) = g_{3\hat{m}}^0(X_{-3\hat{m}}^0(\vec{p}^{\max}))$.*

PROOF.    We first prove the statement for $\mathcal{H}_n(\mathcal{S}, \varepsilon)$. As mentioned before, for $\varepsilon$ sufficiently small ($\varepsilon < \varepsilon_0$), $\mathcal{H}_n(\mathcal{S}, \varepsilon)$ is achieved by $X_{-n}^0$ with $\mathcal{A}(X_{-n}^0) = \mathcal{S}_{n+1}$; and one may choose $\delta$ such that

$$\mathcal{H}_n(\mathcal{S}, \varepsilon) = \sup_{\vec{p}: X_{-n}^0(\vec{p}_n) \in \mathbb{P}_{n+1,\delta}, \mathcal{A}(X_{-n}^0(\vec{p}_n)) = \mathcal{S}_{n+1}} H(Z_0(\vec{p}_n, \varepsilon)|Z_{-n}^{-1}(\vec{p}_n, \varepsilon)).$$

Below, we assume $\varepsilon < \varepsilon_0$, $X_{-n}^0(\vec{p}_n) \in \mathbb{P}_{n+1,\delta}$, $\mathcal{A}(X_{-n}^0(\vec{p}_n)) = \mathcal{S}_{n+1}$; and for notational convenience, we rewrite $f_n^0(X_{-n}^0(\vec{p}_n))$ as $f_n(\vec{p}_n)$, $g_n^0(X_{-n}^0(\vec{p}_n))$ as $g_n(\vec{p}_n)$.

In Lemma 2.2 we have proved that

$$H(Z_0(\vec{p}_n, \varepsilon)|Z_{-n}^{-1}(\vec{p}_n, \varepsilon))$$
$$= H(X_0(\vec{p}_n)|X_{-n}^{-1}(\vec{p}_n)) + f_n(\vec{p}_n)\varepsilon \log(1/\varepsilon) + g_n(\vec{p}_n)\varepsilon + O(\varepsilon^2 \log \varepsilon).$$

Moreover, by Remark 2.3, for any $\delta > 0$, $O(\varepsilon^2 \log \varepsilon)$ is uniform on $\mathbb{P}_{n+1,\delta}$, that is, there is a constant $C$ (depending on $n$) such that, for all $X_{-n}^0$ with

$X^0_{-n}(\vec{p}) \in \mathbb{P}_{n+1,\delta}$ and $\mathcal{A}(X^0_{-n}) = \mathcal{S}_{n+1}$,

$$|H(Z_0(\vec{p}_n, \varepsilon)|Z^{-1}_{-n}(\vec{p}_n, \varepsilon))$$
$$- H(X_0(\vec{p}_n)|X^{-1}_{-n}(\vec{p}_n)) - f_n(\vec{p}_n)\varepsilon \log(1/\varepsilon) - g_n(\vec{p}_n)\varepsilon|$$
$$\leq C\varepsilon^2 \log \varepsilon.$$

Let $\vec{q}_n = \vec{p}_n - \vec{p}^{\max}_n$. Since $H(X_0(\vec{p}_n)|X^{-1}_{-n}(\vec{p}_n))$ is maximized at $\vec{p}^{\max}_n$, we can expand $H(X_0(\vec{p}_n)|X^{-1}_{-n}(\vec{p}_n))$ around $\vec{p}^{\max}_n$:

$$H(X_0(\vec{p}_n)|X^{-1}_{-n}(\vec{p}_n)) = H(X_0(\vec{p}^{\max}_n)|X^{-1}_{-n}(\vec{p}^{\max}_n)) + \vec{q}_n^{\,t} K_1 \vec{q}_n + O(|\vec{q}_n|^3)$$
$$= H(X_{\vec{p}^{\max}}) + \vec{q}_n^t K_1 \vec{q}_n + O(|\vec{q}_n|^3),$$

where $K_1$ is a negative definite matrix by Lemma 3.1 (the second equality follows from the fact that $X_{\vec{p}^{\max}}$ is an $\hat{m}$th order Markov chain). So for $|\vec{q}_n|$ sufficiently small, we have

$$H(X_0(\vec{p}_n)|X^{-1}_{-n}(\vec{p}_n)) < H(X_{\vec{p}^{\max}}) + (1/2)\vec{q}_n^t K_1 \vec{q}_n.$$

Now we expand $f_n(\vec{p}_n)$ and $g_n(\vec{p}_n)$ around $\vec{p}^{\max}_n$:

$$f_n(\vec{p}_n) = f_n(\vec{p}^{\max}_n) + K_2 \cdot \vec{q}_n + O(|\vec{q}_n|^2),$$
$$g_n(\vec{p}_n) = g_n(\vec{p}^{\max}_n) + K_3 \cdot \vec{q}_n + O(|\vec{q}_n|^2)$$

(here, $K_2$ and $K_3$ are vectors of first order partial derivatives). Then, for $|\vec{q}_n|$ sufficiently small, we have

$$f_n(\vec{p}_n) \leq f_n(\vec{p}^{\max}_n) + 2\sum_j |K_{2,j}||q_{n,j}|,$$
$$g_n(\vec{p}) \leq g_n(\vec{p}^{\max}_n) + 2\sum_j |K_{3,j}||q_{n,j}|,$$

where $K_{2,j}, K_{3,j}, q_{n,j}$ are the $j$th coordinates of $K_2, K_3, \vec{q}_n$, respectively.

With a change of coordinates, if necessary, we may assume $K_1$ is a diagonal matrix with strictly negative diagonal elements $K_{1,j}$. In the following we assume $0 < \varepsilon < \varepsilon_0$. And we may further assume that for some $\ell \geq 1$, $|q_{n,j}| > 4|K_{2,j}/K_{1,j}|\varepsilon \log(1/\varepsilon) + 4|K_{3,j}/K_{1,j}|\varepsilon$ for $j \leq \ell - 1$, and $|q_{n,j}| \leq 4|K_{2,j}/K_{1,j}|\varepsilon \log(1/\varepsilon) + 4|K_{3,j}/K_{1,j}|\varepsilon$ for $j \geq \ell$. Then for each $j \leq l - 1$, we have $(1/2)K_{1,j}q_{n,j}^2 + 2|K_{2,j}||q_{n,j}|\varepsilon \log(1/\varepsilon) + 2|K_{3,j}||q_{n,j}|\varepsilon < 0$. Thus,

$$H(Z_0(\vec{p}_n, \varepsilon)|Z^{-1}_{-n}(\vec{p}_n, \varepsilon))$$
$$< H(X_{\vec{p}^{\max}_n}) + f_n(\vec{p}^{\max}_n)\varepsilon \log(1/\varepsilon) + g_n(\vec{p}^{\max}_n)\varepsilon$$
$$+ \sum_j ((1/2)K_{1,j}q_{n,j}^2 + 2|K_{2,j}||q_{n,j}|\varepsilon \log(1/\varepsilon) + 2|K_{3,j}||q_{n,j}|\varepsilon)$$

$$+ C\varepsilon^2 \log \varepsilon$$

$$< H(X_{\vec{p}^{\max}}) + f_n(\vec{p}_n^{\max})\varepsilon \log(1/\varepsilon)$$

$$+ g_n(\vec{p}_n^{\max})\varepsilon + \sum_{j \geq l}(1/2)K_{1,j}(4|K_{2,j}/K_{1,j}|\varepsilon \log(1/\varepsilon) + 4|K_{3,j}/K_{1,j}|\varepsilon)^2$$

$$+ \sum_{j \geq l} 2|K_{2,j}|(4|K_{2,j}/K_{1,j}|\varepsilon \log(1/\varepsilon) + 4|K_{3,j}/K_{1,j}|\varepsilon)\varepsilon \log(1/\varepsilon)$$

$$+ \sum_{j \geq l} 2|K_{3,j}|(4|K_{2,j}/K_{1,j}|\varepsilon \log(1/\varepsilon) + 4|K_{3,j}/K_{1,j}|\varepsilon)\varepsilon + C\varepsilon^2 \log \varepsilon.$$

Collecting terms, we eventually reach

$$H(Z_0(\vec{p}_n, \varepsilon)|Z_{-n}^{-1}(\vec{p}_n, \varepsilon))$$

$$< H(X_{\vec{p}^{\max}}) + f_n(\vec{p}_n^{\max})\varepsilon \log(1/\varepsilon) + g_n(\vec{p}_n^{\max})\varepsilon + O(\varepsilon^2 \log^2 \varepsilon),$$

and since $\mathcal{H}_n(\mathcal{S}, \varepsilon)$ is the sup of the left-hand side expression, together with $H(X_{\vec{p}^{\max}}) = C(\mathcal{S})$, we have

$$\mathcal{H}_n(\mathcal{S}, \varepsilon) \leq C(\mathcal{S}) + f_n(\vec{p}_n^{\max})\varepsilon \log(1/\varepsilon) + g_n(\vec{p}_n^{\max})\varepsilon + O(\varepsilon^2 \log^2 \varepsilon).$$

As discussed in Theorem 2.4, we have

$$(3.3) \qquad f_n(\vec{p}_n^{\max}) = f(X_{\vec{p}^{\max}}), \qquad n \geq 2\hat{m},$$

and

$$(3.4) \qquad g_n(\vec{p}_n^{\max}) = g(X_{\vec{p}^{\max}}), \qquad n \geq 3\hat{m}.$$

So eventually we reach

$$\mathcal{H}_n(\mathcal{S}, \varepsilon) \leq C(\mathcal{S}) + f(X_{\vec{p}^{\max}})\varepsilon \log(1/\varepsilon) + g(X_{\vec{p}^{\max}})\varepsilon + O(\varepsilon^2 \log^2 \varepsilon).$$

The reverse inequality follows trivially from the definition of $\mathcal{H}_n(\varepsilon)$.

We now prove the statement for $h_m(\mathcal{S}, \varepsilon)$. First, observe that

$$\mathcal{H}_{3m}(\mathcal{S}, \varepsilon) \geq h_m(\mathcal{S}, \varepsilon) \geq h_{\hat{m}}(\mathcal{S}, \varepsilon) \geq H(Z_{\vec{p}^{\max}, \varepsilon}),$$

where $Z_{\vec{p}^{\max}, \varepsilon}$ is the output process corresponding to input process $X_{\vec{p}^{\max}}$. By part 1, $\mathcal{H}_{3m}(\mathcal{S}, \varepsilon)$ is of the form $C(\mathcal{S}) + f(X_{\vec{p}^{\max}})\varepsilon \log(1/\varepsilon) + g(X_{\vec{p}^{\max}})\varepsilon + O(\varepsilon^2 \log^2 \varepsilon)$. By Theorem 2.4, $H(Z_{\vec{p}^{\max}, \varepsilon})$ is of the same form. Thus, $h_m(\mathcal{S}, \varepsilon)$ is also of the same form, as desired. $\square$

**COROLLARY 3.3.**

$$C(\mathcal{S}, \varepsilon) = C(\mathcal{S}) + (f(X_{\vec{p}^{\max}}) - 1)\varepsilon \log(1/\varepsilon) + (g(X_{\vec{p}^{\max}}) - 1)\varepsilon + O(\varepsilon^2 \log^2 \varepsilon).$$

*In fact, for each $m \geq \hat{m}(\mathcal{S})$, $h_m(S, \varepsilon) - H(\varepsilon)$ is of this form.*

PROOF. This follows from Theorem 3.2, inequality (3.1) and the fact that

$$H(\varepsilon) = \varepsilon \log 1/\varepsilon + (1-\varepsilon)\log 1/(1-\varepsilon) = \varepsilon \log 1/\varepsilon + \varepsilon + O(\varepsilon^2). \qquad \square$$

REMARK 3.4. Note that the error term here for noisy constrained capacity is $O(\varepsilon^2 \log^2 \varepsilon)$, which is larger than the error term, $O(\varepsilon^2 \log \varepsilon)$, for entropy rate in Theorem 2.4. At least in some cases, this cannot be improved, as we show at the end of the next section.

**4. Binary symmetric channel with $(d,k)$-RLL constrained input.** We now apply the results of the preceding section to compute asypmtotics for the the noisy constrained BSC channel with inputs restricted to the $(d,k)$-RLL constraint $\mathcal{S}(d,k)$. Expressions (2.8) and (2.9) allow us to explicitly compute $f(X_{\overline{p}^{\max}})$ and $g(X_{\overline{p}^{\max}})$. In this section, as an example, we derive the explicit expression for $f(X_{\overline{p}^{\max}})$, omitting the computation of $g(X_{\overline{p}^{\max}})$ due to tedious derivation. We remark that for a BSC($\varepsilon$) for some cases, the $(d,k)$-RLL constrained input, similar expressions have been independently obtained in [19].

It is first shown in [19] that in the case $k \leq 2d$, for any binary stationary Markov chain $X$, of any order, with $\mathcal{A}(X) \subseteq \mathcal{S}(d,k)$, $f(X) = 1$, and so, in this case, $C(\mathcal{S}(d,k), \varepsilon) = C(\mathcal{S}(d,k), 0) + O(\varepsilon)$, that is, the noisy constrained capacity differs from the noiseless capacity by $O(\varepsilon)$, rather than $O(\varepsilon \log \varepsilon)$. In the following we take a look at this using a different approach. For this, first note that for any $d, k$, $f(X)$ takes the form

$$\begin{aligned}
(4.1) \quad f(X) &= \sum_{l_1+l_2 \leq k-1, 0 \leq l_2 \leq d-1, l_1 \geq d} p_X(10^{l_1+l_2+1}1) \\
&\quad + \sum_{l_1+l_2=k, l_1 \geq d} p_X(10^{l_1}10^{l_2}) + \sum_{1 \leq l \leq k} p_X(10^l).
\end{aligned}$$

Now, when $k \leq 2d$,

$$\sum_{l_1+l_2=k, l_1 \geq d} p_X(10^{l_1}10^{l_2}) = \sum_{d \leq l_1 \leq k} p_X(10^{l_1}1) = p(1)$$

and

$$\begin{aligned}
\sum_{l_1+l_2 \leq k-1, 0 \leq l_2 \leq d-1, l_1 \geq d} &p_X(10^{l_1+l_2+1}1) \\
&= p_X(10^{d+1}) + p_X(10^{d+2}) + \cdots + p_X(10^k).
\end{aligned}$$

So

$$f(X) = p_X(1) + p_X(10) + \cdots + p_X(10^d) + p_X(10^{d+1}) + \cdots + p_X(10^k) = 1,$$

as desired.

Now we consider the general RLL constraint $\mathcal{S}(d,k)$. By Corollary 3.3, we have

$$
\begin{aligned}
(4.2) \quad C(\mathcal{S}(d,k),\varepsilon) = {}& C(\mathcal{S}(d,k)) + (f(X_{\vec{p}_{\max}}) - 1)\varepsilon \log 1/\varepsilon \\
& + (g(X_{\vec{p}_{\max}}) - 1)\varepsilon + O(\varepsilon^2 \log^2 \varepsilon).
\end{aligned}
$$

For any irreducible finite type constraint, the noiseless capacity and Markov process of maximal entropy rate can be computed in various ways (which all go back to Shannon; see [21] or [20], page 444). Let $A$ denote the adjacency matrix of the standard graph presentation, with $k+1$ states, of $\mathcal{S}(d,k)$. Let $\rho$ denote the reciprocal of the largest eigenvalue. One can write $C(\mathcal{S}(d,k)) = -\log \rho_0$, and in this case $\rho_0$ is the real root of

$$
(4.3) \qquad\qquad \sum_{\ell=d}^{k} \rho_0^{\ell+1} = 1.
$$

In the following we compute $f(X_{\vec{p}_{\max}})$ explicitly in terms of $\rho_0$. Let $\vec{w} = (w_0, w_1, \ldots, w_k)$ and $\vec{v} = (v_0, v_1, \ldots, v_k)$ denote the left and right eigenvectors of $A$. Assume that $\vec{w}$ and $\vec{v}$ are scaled such that $\vec{w} \cdot \vec{v} = 1$. Then one checks that, with $X = X_{\vec{p}^{\max}}$,

$$
p_X(1) = w_0 v_0 = \frac{1}{(k+1) - \sum_{j=d+1}^{k} \sum_{l=0}^{j-d-1} 1/\rho_0^{l-j}},
$$

$$
p_X(10^{l_1 + l_2 + 1} 1) = p_X(1)\rho_0^{l_1 + l_2 + 2}, \qquad p_X(10^k 1) = p_X(1)\rho_0^{k+1},
$$

$$
\begin{aligned}
p_X(10^{l_1} 10^{l_2}) &= p_X(10^{l_1} 10^{l_2} 1) + p_X(10^{l_1} 10^{l_2+1} 1) + \cdots + p_X(10^{l_1} 10^k 1) \\
&= p_X(1)\rho_0^{l_1 + l_2 + 2}(1 + \rho_0 + \cdots + \rho_0^{k-l_2}) \\
&= p_X(1)\rho_0^{l_1 + l_2 + 2}\frac{1 - \rho_0^{k-l_2+1}}{1 - \rho_0}
\end{aligned}
$$

and

$$
\begin{aligned}
p_X(10^l) &= p_X(10^l 1) + p_X(10^{l+1} 1) + \cdots + p_X(10^k 1) \\
&= p_X(1)\rho_0^{l+1}(1 + \rho_0 + \cdots + \rho_0^{k-l}) = p_X(1)\rho_0^{l+1}\frac{1 - \rho_0^{k-l+1}}{1 - \rho_0}.
\end{aligned}
$$

So we obtain an explicit expression:

$$
\begin{aligned}
f(X_{\vec{p}^{\max}}) = {}& \sum_{l_1 + l_2 \le k-1, 0 \le l_2 \le d-1, l_1 \ge d} p_X(10^{l_1 + l_2 + 1} 1) \\
& + \left( \sum_{l_1 = k, l_2 = 0} + \sum_{l_1 + l_2 = k, k-1 \ge l_1 \ge d} \right) p_X(10^{l_1} 10^{l_2}) + \sum_{1 \le l \le d} p_X(10^l)
\end{aligned}
$$

$$= p_X(1)\rho_0^{k+1} + \sum_{l_1+l_2 \leq k-1, 0 \leq l_2 \leq d-1, l_1 \geq d} p_X(1)\rho_0^{l_1+l_2+2}$$

$$+ \sum_{l_1+l_2=k, k-1 \geq l_1 \geq d} p_X(1)\rho_0^{l_1+l_2+2}\frac{1-\rho_0^{k-l_2+1}}{1-\rho_0}$$

$$+ \sum_{1 \leq l \leq d} p_X(1)\rho_0^{l+1}\frac{1-\rho_0^{k-l+1}}{1-\rho_0}.$$

The coefficient $g$ can also be computed explicitly but takes a more complicated form.

EXAMPLE 4.1. Consider a first order stationary Markov chain $X$ with $\mathcal{A}(X) \subseteq \mathcal{S}(1, \infty)$, transmitted over BSC($\varepsilon$) with the corresponding output $Z$, a hidden Markov chain. In this case, $X$ can be characterized by the following probability vector:

$$\vec{p}_1 = (p_X(00), p_X(01), p_X(10)).$$

Note that $\hat{m}(\mathcal{S}) = 1$, and the only sequence $w_{-2}w_{-1}v$, which satisfies the requirement that $w_{-2}w_{-1}$ is in $\mathcal{S}$ and $w_{-2}w_{-1}v$ is not allowable in $\mathcal{S}$, is 011. It then follows that

$$(4.4) \quad f(X_{\vec{p}_1}) = p(01\bar{1}) + p(0\bar{1}1) + p(\bar{0}11) = \pi_{01}(2-\pi_{01})/(1+\pi_{01}),$$

where $\pi_{01}$ denotes the transition probability from 0 to 1 in $X$. Straightforward, but tedious, computation also leads to

$$g(X_{\vec{p}_1}) = (1+\pi_{01})^{-1}(2\pi_{01} - \pi_{01}^2 - 2\pi_{01}^3 + 3\pi_{01}^4 - \pi_{01}^5$$
$$+ (-2\pi_{01} + 4\pi_{01}^3 - 2\pi_{01}^4)\log(2)$$
$$+ (-1 + 3\pi_{01} - \pi_{01}^2 - 2\pi_{01}^3 + 5\pi_{01}^4 - 3\pi_{01}^5)\log(\pi_{01})$$
$$+ (2 - 6\pi_{01} + 7\pi_{01}^3 - 8\pi_{01}^4 + 3\pi_{01}^5)\log(1-\pi_{01})$$
$$+ (2\pi_{01} + \pi_{01}^2 - 3\pi_{01}^3 + \pi_{01}^4)\log(2-\pi_{01})).$$

Thus,

$$H(Z_{\vec{p}_1, \varepsilon}) = H(X_{\vec{p}_1}) + (\pi_{01}(2-\pi_{01})/(1+\pi_{01}))\varepsilon \log(1/\varepsilon)$$
$$+ (g(X_{\vec{p}_1}) - 1)\varepsilon + O(\varepsilon^2 \log \varepsilon).$$

This asymptotic formula was originally proven in [23], with the less precise result that replaces $(g(X_{\vec{p}_1}) - 1)\varepsilon + O(\varepsilon^2 \log^2(1/\varepsilon))$ by $O(\varepsilon)$.

The maximum entropy Markov chain $X_{\vec{p}^{\max}}$ on $\mathcal{S}(1, \infty)$ is defined by the transition probability matrix

$$\begin{bmatrix} 1/\lambda & 1/\lambda^2 \\ 1 & 0 \end{bmatrix}$$

and

$$C(\mathcal{S}) = H(X_{\vec{p}^{\max}}) = \log \lambda,$$

where $\lambda$ is the golden mean. Thus, in this case $\pi_{01} = 1/\lambda^2$ and so by Corollary 3.3, we obtain

$$C(\mathcal{S}, \varepsilon) = \log \lambda - ((2\lambda + 2)/(4\lambda + 3))\varepsilon \log(1/\varepsilon)$$
$$+ (g(X_{\vec{p}_1})|_{\pi_{01}=1/\lambda^2} - 1)\varepsilon + O(\varepsilon^2 \log^2(1/\varepsilon)).$$

We now show that the error term in the above formula cannot be improved, in the sense that the error term is of size at least a positive constant times $\varepsilon^2 \log^2(1/\varepsilon)$. First observe that if we parameterize $\vec{p}_1 = \vec{p}_1(\varepsilon)$ in any way, we obtain

$$(4.5) \qquad C(\mathcal{S}, \varepsilon) \geq H(Z_{\vec{p}_1(\varepsilon),\varepsilon}) - H(\varepsilon).$$

Since $\vec{p}_1$ is uniquely determined by the transition probability $\pi_{01}$, we shall re-write $\vec{p}_1(\varepsilon)$ as $\pi_{01}(\varepsilon)$. We shall also rewrite the value of $\pi_{01} = 1/\lambda^2$ at the maximum entropy Markov chain as $p_{\max}$.

Choose the parametrization $\pi_{01}(\varepsilon) = p_{\max} + \alpha\varepsilon \log(1/\varepsilon)$, where $\alpha$ is selected as follows. Let $K_1$ denote the value of the second derivative of $H(X_{\pi_{01}})$ at $\pi_{01} = p_{\max}$ (the first derivative at $\pi_{01} = p_{\max}$ is 0). Let $K_2$ denote the value of the first derivative of $f(X_{\pi_{01}})$ at $\pi_{01} = p_{\max}$. These values can be computed explicitly: $K_1$ from the formula for entropy rate of a first order Markov chain (1.3) and $K_2$ from (4.4) above. A computation shows that $K_1 \approx -3.065$ and $K_2 \approx 0.571$ (all that really matters is that neither constant is 0). Let $\alpha$ be any number such that $0 < \alpha < K_2/|K_1|$.

From Theorem 2.4 and Remark 2.6, we have

$$(4.6) \qquad H(Z_{\pi_{01}(\varepsilon),\varepsilon}) \geq H(X_{\pi_{01}(\varepsilon)}) + f(X_{\pi_{01}(\varepsilon)})\varepsilon \log(1/\varepsilon)$$
$$+ g(X_{\pi_{01}(\varepsilon)})\varepsilon + C_1\varepsilon^2 \log \varepsilon,$$

for some constant $C_1$ (independent of $\varepsilon$ sufficiently small). We also have

$$(4.7) \quad H(X_{\pi_{01}(\varepsilon)}) \geq H(X_{p_{\max}}) + K_1(\alpha\varepsilon \log(1/\varepsilon))^2 + C_2(\alpha\varepsilon \log(1/\varepsilon))^3$$

for some constant $C_2$. And

$$(4.8) \quad f(X_{\pi_{01}(\varepsilon)}) \geq f(X_{p_{\max}}) + K_2(\alpha\varepsilon \log(1/\varepsilon)) + C_3(\alpha\varepsilon \log(1/\varepsilon))^2,$$

$$(4.9) \quad g(X_{\pi_{01}(\varepsilon)}) \geq g(X_{p_{\max}}) + C_4(\alpha\varepsilon \log(1/\varepsilon))$$

for constants $C_3, C_4$. And recall that

$$(4.10) \quad H(\varepsilon) = \varepsilon \log 1/\varepsilon + (1-\varepsilon)\log 1/(1-\varepsilon) = \varepsilon \log 1/\varepsilon + \varepsilon + C_5\varepsilon^2$$

for some constant $C_5$.

Recalling that $H(X_{p_{\max}}) = C(\mathcal{S})$ and combining (4.5)–(4.10), we see that

$$C(\mathcal{S}, \varepsilon) \geq C(\mathcal{S}) + (f(X_{p_{\max}}) - 1)\varepsilon \log(1/\varepsilon) + (g(X_{p_{\max}}) - 1)\varepsilon$$
$$+ K_1(\alpha\varepsilon \log(1/\varepsilon))^2 + K_2(\alpha\varepsilon^2 \log^2(1/\varepsilon))$$

plus "error terms" which add up to

$$C_1\varepsilon^2 \log\varepsilon + C_2(\alpha\varepsilon \log(1/\varepsilon))^3 + C_3\alpha^2(\varepsilon \log(1/\varepsilon))^3$$
$$+ C_4(\alpha\varepsilon^2 \log(1/\varepsilon)) + C_5\varepsilon^2,$$

which is lower bounded by a constant $M$ times $\varepsilon^2 \log(1/\varepsilon)$. Thus, we see that the difference between $C(\mathcal{S}, \varepsilon)$ and $C(S) + (f(X_{p_{\max}}) - 1)\varepsilon \log(1/\varepsilon) + (g(X_{p_{\max}}) - 1)\varepsilon$ is lower bounded by

$$(4.11) \qquad \alpha(K_1\alpha + K_2)\varepsilon^2 \log^2(1/\varepsilon) + M\varepsilon^2 \log(1/\varepsilon).$$

Since $\alpha > 0$ and $K_1\alpha + K_2 > 0$, for sufficiently small $\varepsilon$, (4.11) is lower bounded by a positive constant times $\varepsilon^2 \log^2(1/\varepsilon)$, as desired.

## APPENDIX

We first prove (2.11)–(2.13).

- (2.11) follows trivially from the fact that $X$ is an $m$th order Markov chain.
- Now consider (2.12). For $w \in \mathcal{A}(X)$ and $wv \notin \mathcal{A}(X)$,

$$h_{n+1}^k(wv) = \sum_{j=1}^{n-k} p_X(w_{-n}^{-j-1}\bar{w}_{-j}w_{-j+1}^{-1}v) + p_X(w_{-n}^{-1}\bar{v})$$

$$= \sum_{j=1}^{m} p_X(w_{-n}^{-j-1}\bar{w}_{-j}w_{-j+1}^{-1}v) + p_X(w_{-n}^{-1}\bar{v}).$$

So

$$\frac{h_{n+1}^k(wv)}{p_X(w)} = \frac{\sum_{j=1}^{m} p_X(w_{-n}^{-j-1}\bar{w}_{-j}w_{-j+1}^{-1}v) + p_X(w_{-n}^{-1}\bar{v})}{p_X(w_{-n}^{-1})}$$

$$= \left( \sum_{j=1}^{m} p_X(w_{-m}^{-j-1}\bar{w}_{-j}w_{-j+1}^{-1}v | w_{-2m}^{-m-1}) \right.$$

$$\left. + p_X(w_{-m}^{-1}\bar{v} | w_{-2m}^{-m-1}) \right)$$

$$\times p_X(w_{-n}^{-m-1})(p_X(w_{-m}^{-1}|w_{-2m}^{-m-1})p_X(w_{-n}^{-m-1}))^{-1}$$

$$= \frac{\sum_{j=1}^{m} p_X(w_{-2m}^{-j-1}\bar{w}_{-j}w_{-j+1}^{-1}v) + p_X(w_{-2m}^{-1}\bar{v})}{p_X(w_{-2m}^{-1})}.$$

- For (2.13), there are two cases. If $p_X(w_{-n}^{-m-1}) = 0$,

$$\frac{h_{n+1}^k(wv)}{h_n^k(w)} = \frac{\sum_{j=1}^{n-k} p_X(w_{-n}^{-j-1} \bar{w}_{-j} w_{-j+1}^{-1} v)}{\sum_{j=1}^{n-k} p_X(w_{-n}^{-j-1} \bar{w}_{-j} w_{-j+1}^{-1})}$$

$$= \frac{\sum_{j=m+1}^{n-k} p_X(w_{-n}^{-j-1} \bar{w}_{-j} w_{-j+1}^{-1} v)}{\sum_{j=m+1}^{n-k} p_X(w_{-n}^{-j-1} \bar{w}_{-j} w_{-j+1}^{-1})} = p_X(v|w_{-m}^{-1}).$$

If $p_X(w_{-n}^{-m-1}) > 0$,

$$\frac{h_{n+1}^k(wv)}{h_n^k(w)} = \frac{\sum_{j=1}^{n-k} p_X(w_{-n}^{-j-1} \bar{w}_{-j} w_{-j+1}^{-1} v)}{\sum_{j=1}^{n-k} p_X(w_{-n}^{-j-1} \bar{w}_{-j} w_{-j+1}^{-1})}$$

$$= \frac{\sum_{j=1}^{2m} p_X(w_{-n}^{-j-1} \bar{w}_{-j} w_{-j+1}^{-1} v)}{\sum_{j=1}^{2m} p_X(w_{-n}^{-j-1} \bar{w}_{-j} w_{-j+1}^{-1})}$$

$$= \frac{\sum_{j=1}^{2m} p_X(w_{-3m}^{-j-1} \bar{w}_{-j} w_{-j+1}^{-1} v)}{\sum_{j=1}^{2m} p_X(w_{-3m}^{-j-1} \bar{w}_{-j} w_{-j+1}^{-1})}.$$

Using (2.11)–(2.13), we now proceed to prove (2.14)–(2.16).

- For (2.14), we have

$$\sum_{wv \in \mathcal{A}(X)} \frac{h_{n+1}^k(wv) p_X(w) - h_n^k(w) p_X(wv)}{p_X(w)}$$

$$= \sum_{wv \in \mathcal{A}(X)} h_{n+1}^k(wv) - \sum_{wv \in \mathcal{A}(X)} h_n^k(w) p_X(v|w_{-m}^{-1})$$

$$= \sum_{wv \in \mathcal{A}(X)} \left( \sum_{j=1}^{n-k} p_X(w_{-n}^{-j-1} \bar{w}_{-j} w_{-j+1}^{-1} v) + p_X(w_{-n}^{-1} \bar{v}) \right)$$

$$\quad - (n+1-k) \sum_{wv \in \mathcal{A}(X)} p_X(wv)$$

$$\quad - \sum_{wv \in \mathcal{A}(X)} \sum_{j=1}^{n-k} p_X(w_{-n}^{-j-1} \bar{w}_{-j} w_{-j+1}^{-1}) p_X(v|w_{-m}^{-1})$$

$$\quad + (n-k) \sum_{wv \in \mathcal{A}(X)} p_X(wv)$$

$$= \sum_{wv \in \mathcal{A}(X)} \sum_{j=1}^{n-k} p_X(w_{-n}^{-j-1} \bar{w}_{-j} w_{-j+1}^{-1} v)$$

$$- \sum_{w \in \mathcal{A}(X)} \sum_{j=1}^{n-k} p_X(w_{-n}^{-j-1} \bar{w}_{-j} w_{-j+1}^{-1}) + \sum_{wv \in \mathcal{A}(X)} p_X(w_{-n}^{-1} \bar{v}) - 1$$

$$= \sum_{wv \in \mathcal{A}(X)} \sum_{j=1}^{n-k} p_X(w_{-n}^{-j-1} \bar{w}_{-j} w_{-j+1}^{-1} v)$$

$$- \sum_{w \in \mathcal{A}(X)} \left( \sum_{j=1}^{n-k} p_X(w_{-n}^{-j-1} \bar{w}_{-j} w_{-j+1}^{-1} 0) \right.$$
$$\left. + \sum_{j=1}^{n-k} p_X(w_{-n}^{-j-1} \bar{w}_{-j} w_{-j+1}^{-1} 1) \right)$$

$$+ \sum_{w_{-m}^{-1} v \in \mathcal{A}(X)} p_X(w_{-m}^{-1} \bar{v}) - 1$$

$$= \sum_{wv \in \mathcal{A}(X)} \sum_{j=1}^{n-k} p_X(w_{-n}^{-j-1} \bar{w}_{-j} w_{-j+1}^{-1} v)$$

$$- \sum_{wv \in \mathcal{A}(X)} \sum_{j=1}^{n-k} p_X(w_{-n}^{-j-1} \bar{w}_{-j} w_{-j+1}^{-1} v)$$

$$- \sum_{w \in \mathcal{A}(X), wv \notin \mathcal{A}(X)} \sum_{j=1}^{n-k} p_X(w_{-n}^{-j-1} \bar{w}_{-j} w_{-j+1}^{-1} v)$$
$$+ \sum_{w_{-m}^{-1} v \in \mathcal{A}(X)} p_X(w_{-m}^{-1} \bar{v}) - 1$$

$$= - \sum_{w_{-2m}^{-1} \in \mathcal{A}(X), w_{-2m}^{-1} v \notin \mathcal{A}(X)} \sum_{j=1}^{m} p_X(w_{-2m}^{-j-1} \bar{w}_{-j} w_{-j+1}^{-1} v)$$
$$+ \sum_{w_{-m}^{-1} v \in \mathcal{A}(X)} p_X(w_{-m}^{-1} \bar{v}) - 1.$$

- For (2.15), we have

$$\sum_{w \in \mathcal{A}(X), wv \notin \mathcal{A}(X)} h_{n+1}^k(wv) \log \frac{h_{n+1}^k(wv)}{p_X(w)}$$
$$= \sum_{w \in \mathcal{A}(X), wv \notin \mathcal{A}(X)} h_{n+1}^k(wv) \log \frac{h_{2m+1}^0(w_{-2m}^{-1} v)}{p_X(w_{-2m}^{-1})}$$

$$
\begin{aligned}
&= \sum_{w \in \mathcal{A}(X), wv \notin \mathcal{A}(X)} \sum_{j=1}^{n-k} p_X(w_{-n}^{-j-1} \bar{w}_{-j} w_{-j+1}^{-1} v) \log \frac{h_{2m+1}^0(w_{-2m}^{-1} v)}{p_X(w_{-2m}^{-1})} \\
&= \sum_{w_{-2m}^{-1} \in \mathcal{A}(X), w_{-2m}^{-1} v \notin \mathcal{A}(X)} \sum_{j=1}^{m} p_X(w_{-2m}^{-j-1} \bar{w}_{-j} w_{-j+1}^{-1} v) \\
&\qquad\qquad\qquad \times \log \frac{h_{2m+1}^0(w_{-2m}^{-1} v)}{p_X(w_{-2m}^{-1})}.
\end{aligned}
$$

- For (2.16), we have

$$
\sum_{wv \in \mathcal{A}(X)} h_{n+1}^k(wv) \log p_X(v|w)
$$

$$
+ \sum_{p_{XZ}(w)=\Theta(\varepsilon), p_{XZ}(wv)=\Theta(\varepsilon)} h_{n+1}^k(wv) \log \frac{h_{n+1}^0(wv)}{h_n^0(w)}
$$

$$
= \sum_{wv \in \mathcal{A}(X)} \left( \sum_{j=1}^{n-k} p_X(w_{-n}^{-j-1} \bar{w}_{-j} w_{-j+1}^{-1} v) \right.
$$

$$
\left. + p_X(w_{-n}^{-1} \bar{v}) - (n+1-k) p_X(wv) \right) \log p_X(v|w_{-m}^{-1})
$$

$$
+ \sum_{p_{XZ}(w)=\Theta(\varepsilon), p_{XZ}(wv)=\Theta(\varepsilon), p_X(w_{-n}^{-m-1})=0} h_{n+1}^k(wv) \log \frac{h_{n+1}^k(wv)}{h_n^k(w)}
$$

$$
+ \sum_{p_{XZ}(w)=\Theta(\varepsilon), p_{XZ}(wv)=\Theta(\varepsilon), p_X(w_{-n}^{-m-1})>0} h_{n+1}^k(wv) \log \frac{h_{n+1}^k(wv)}{h_n^k(w)}
$$

$$
= \left( \sum_{wv \in \mathcal{A}(X)} + \sum_{p_{XZ}(w)=\Theta(\varepsilon), p_{XZ}(wv)=\Theta(\varepsilon), p_X(w_{-n}^{-m-1})=0} \right)
$$

$$
\left( \sum_{j=1}^{n-k} p_X(w_{-n}^{-j-1} \bar{w}_{-j} w_{-j+1}^{-1} v) + p_X(w_{-n}^{-1} \bar{v}) \right) \log p_X(v|w_{-m}^{-1})
$$

$$
- (n+1-k) \sum_{w_{-m}^{-1} v \in \mathcal{A}(X)} p_X(w_{-m}^{-1} v) \log p_X(v|w_{-m}^{-1})
$$

$$
+ \sum_{p_{XZ}(w_{-3m}^{-1})=\Theta(\varepsilon), p_{XZ}(w_{-3m}^{-1} v)=\Theta(\varepsilon), p_X(w_{-3m}^{-m-1})>0} h_{3m+1}^0(wv)
$$

$$
\times \log \frac{h_{3m+1}^0(wv)}{h_{3m}^0(w)}
$$

$$= (n - k - m) \sum_{w_{-m}^{-1} v \in \mathcal{A}(X)} p_X(w_{-m}^{-1}v) \log p_X(v|w_{-m}^{-1})$$

$$+ \sum_{wv \in \mathcal{A}(X)} \left( \sum_{j=1}^{m} p(w_{-n}^{-j-1}\bar{w}_{-j}w_{-j+1}^{-1}v) + p_X(w_{-n}^{-1}\bar{v}) \right) \log p_X(v|w_{-m}^{-1})$$

$$- (n + 1 - k) \sum_{w_{-m}^{-1} v \in \mathcal{A}(X)} p_X(w_{-m}^{-1}v) \log p_X(v|w_{-m}^{-1})$$

$$+ \sum_{p_{XZ}(w_{-3m}^{-1})=\Theta(\varepsilon), p_{XZ}(w_{-3m}^{-1}v)=\Theta(\varepsilon), p_X(w_{-3m}^{-m-1})>0} h_{3m+1}^0(wv)$$

$$\times \log \frac{h_{3m+1}^0(wv)}{h_{3m}^0(w)}$$

$$= (-m - 1) \sum_{w_{-m}^{-1} v \in \mathcal{A}(X)} p_X(w_{-m}^{-1}v) \log p_X(v|w_{-m}^{-1})$$

$$+ \sum_{w_{-2m}^{-1} v \in \mathcal{A}(X)} \left( \sum_{j=1}^{m} p_X(w_{-2m}^{-j-1}\bar{w}_{-j}w_{-j+1}^{-1}v) + p_X(w_{-2m}^{-1}\bar{v}) \right)$$

$$\times \log p_X(v|w_{-m}^{-1})$$

$$+ \sum_{p_{XZ}(w_{-3m}^{-1})=\Theta(\varepsilon), p_{XZ}(w_{-3m}^{-1}v)=\Theta(\varepsilon), p_X(w_{-3m}^{-m-1})>0} h_{3m+1}^0(wv)$$

$$\times \log \frac{h_{3m+1}^0(wv)}{h_{3m}^0(w)}.$$

**Acknowledgments.** We are grateful to Wojciech Szpankowski, who raised the problem addressed in this paper and suggested a version of the result in Corollary 3.3. We also thank the anonymous reviewers for helpful comments.

DEPARTMENT OF MATHEMATICS
UNIVERSITY OF HONG KONG
POKFULAM ROAD, POKFULAM
HONG KONG
E-MAIL: ghan@maths.hku.hk

DEPARTMENT OF MATHEMATICS
UNIVERSITY OF BRITISH COLUMBIA
VANCOUVER, B.C. V6T 1Z2
CANADA
E-MAIL: marcus@math.ubc.ca